\documentclass{amsart}
\usepackage{graphicx,amssymb,stmaryrd,amsfonts,epsfig,amsthm,a4,amsmath,mathrsfs,url}
\usepackage{vmargin,xypic,amscd}
\usepackage{eufrak}
\usepackage[latin1]{inputenc}
\psfigdriver{dvips}

\newtheorem{thm}{Theorem}

\newtheorem{cor}[thm]{Corollary}

\newtheorem{prop}[thm]{Proposition}
\theoremstyle{definition}

\theoremstyle{remark}
\newtheorem{rem}[thm]{Remark}
\newtheorem{que}[thm]{Question}

\newcommand{\bpr}{\noindent \textbf{Proof}: ~}
\newcommand{\epr}{~$\blacksquare$}
\begin{document}


\title[Quotients of word hyperbolic groups]
{A note on quotients of word hyperbolic groups with Property~(T)}
\author{Yves de Cornulier}%
\date{\today}
\begin{abstract}
Every discrete group with Kazhdan's Property~(T) is a quotient of
a torsion-free, word hyperbolic group with Property~(T).
\end{abstract}
\maketitle

All groups in the paper are discrete and countable. Recall
\cite{HV} that a group $G$ has Property~(T) if every isometric
action of $G$ on an affine Hilbert space has a fixed point (or,
equivalently, has bounded orbits).

As an immediate consequence of the definition, Property~(T) is
inherited by quotients and by extensions.

Shalom \cite[Theorem 6.7]{Sh} has proved the following interesting
result about Property~(T).

\begin{thm}[Shalom, 2000]
For every group $G$ with Property~(T), there exists a finitely
presented group $G_0$ with Property~(T) which maps onto
$G$.\label{Sha}
\end{thm}

In other words, this means that, given a finite generating subset
for $G$, only finitely many relations suffice to imply
Property~(T). This can be interpreted in the topology of marked
groups \cite{Cha} as: Property~(T) is an open property. See
\cite[3.8]{Gr} for a generalization to other fixed point
properties.

\medskip

A word hyperbolic group is a finitely generated group whose Cayley
graph satisfies a certain condition, introduced by Gromov, meaning
that, at large scale, it is negatively curved. We refer to
\cite{GH} for a precise definition that we do not need here. We
only mention here that word hyperbolic groups are necessarily
finitely presented, that word hyperbolicity is a fundamental
notion in combinatorial group theory as in geometric topology.
Word hyperbolic groups are groups with ``many" quotients, and thus
can be considered as a generalization of free groups.

\medskip

It was asked \cite[Question 16]{W} whether every group with
Property~(T) is quotient of a group with Property~(T) with
finiteness conditions stronger than finite presentation. We give
an answer here by showing that we can impose word hyperbolicity.

\begin{prop}
For every group $G$ with Property~(T), there exists a torsion-free
word hyperbolic group $G_0$ with Property~(T) which maps onto
$G$.\label{FA_T_quot_hyp}
\end{prop}

Note that Proposition \ref{FA_T_quot_hyp} contains Theorem
\ref{Sha} as a corollary; however it is proved by combining
Theorem \ref{Sha} with the following remarkable result of Ollivier
and Wise \cite{OW}. Since it involves some technical definitions,
we do not quote it in full generality.

\begin{thm}[Ollivier and Wise, 2005]
\label{thm:TRips} To every finitely presented group $Q$, we can
associate a short exact sequence $1\rightarrow N \rightarrow G
\rightarrow Q\rightarrow 1$ such that \begin{enumerate}
 \item $G$ is torsion-free, word hyperbolic,

\item $N$ is 2-generated and has property~(T).
\end{enumerate}
\end{thm}

\begin{cor}
For every finitely presented group $Q$ with Property~(T), there
exists a torsion-free word-hyperbolic group $G$ with Property~(T)
mapping onto $Q$ with finitely generated kernel.\label{cor1}
\end{cor}
\bpr Apply Theorem \ref{thm:TRips} to $Q$, so that $G$ lies in an
extension $1\to N\to G\to Q\to 1$, where $N$ has Property~(T) and
$Q$ has Property~(T). Since $Q$ has Property~(T) and since
Property~(T) is stable under extensions, $G$ also has Property
(T).\epr

\begin{rem}
Corollary \ref{cor1} answers a question at the end of \cite{OW}.
\end{rem}

\noindent \textbf{Proof of Proposition \ref{FA_T_quot_hyp}}: Let
$G$ be a group with Property~(T). By Theorem~\ref{Sha}, there
exists a finitely presented group $Q$ with Property~(T) mapping
onto $G$, and by Corollary \ref{cor1}, there exists a torsion-free
word hyperbolic group $G_0$ with Property~(T) mapping onto $Q$, so
that $G_0$ maps onto~$G$.\epr

\begin{que}
1) In Theorem \ref{thm:TRips}, can $G$ be chosen, in addition,
residually finite? In \cite{Wise}, a similar result is proved, $G$
being torsion-free, word hyperbolic, residually finite, and $N$
finitely generated, but never having Property~(T).

2) Let $G$ be a word hyperbolic group (maybe torsion-free), and
$H$ a quotient of $G$ generated by $r$ elements. Does there exist
an intermediate quotient which is both word hyperbolic and
generated by $r$ elements? (The analog statement with ``word
hyperbolic" replaced by ``finitely presented" is immediate.) The
motivation is that, in Proposition \ref{FA_T_quot_hyp}, we would
like to have $G_0$ generated by no more elements than $G$. Theorem
\ref{thm:TRips} only tells us that if $G$ is $r$-generated, then
$G_0$ can be chosen $(r+2)$-generated.
\end{que}

\begin{que}
Following \cite{Se}, a group has Property (FA) if every isometric
action on a simplicial tree has a fixed point. Is it true that
every group with Property (FA) is a quotient of a finitely
presented group with Property (FA)? The author imprudently claimed
that this is true in an earlier version of this paper, but the
proof contained a mistake.
\end{que}


\bigskip

\footnotesize
\noindent Yves de Cornulier\\
\'Ecole Polytechnique Fédérale de Lausanne (EPFL)\\
Institut de Géométrie, Algèbre et Topologie (IGAT)\\
CH-1015 Lausanne, Switzerland\\
E-mail: \url{decornul@clipper.ens.fr}\\

\end{document}